\documentclass[12pt]{article}
\usepackage{cite,amsfonts,amssymb}

\newtheorem{theorem}{Theorem}
\newtheorem{corollary}{Corollary}
\newtheorem{lemma}{Lemma}
\newtheorem{proposition}{Proposition}

\newtheorem{definition}{Definition}
\newtheorem{remark}{Remark}

\title
 {\bf On stochastic fractional Volterra equations in Hilbert space}

\author{{\large\sf Anna Karczewska and Carlos Lizama}\\[2mm]
  \normalsize\it
 Department of Mathematics,
 University of Zielona G\'ora\\ \normalsize\it
 ul. Szafrana 4a, 65-246 Zielona G\'ora, Poland,~ 
 e-mail: A.Karczewska@im.uz.zgora.pl\\[2mm] \normalsize\it
 Universidad de Santiago de Chile, Departamento de
Matem\'atica, Facultad de Ciencias\\ \normalsize\it
Casilla 307-Correo 2,Santiago, Chile,~ 
e-mail: clizama@lauca.usach.cl  }

\begin{document}
\maketitle

\def\thefootnote{}
\footnotetext{
The second author is partially financed by FONDECYT  Grant \#1050084 }

\begin{abstract}
 In this paper, stochastic Volterra equations, particularly fractional,
 in Hilbert space are studied. Sufficient conditions for mild solutions
 to be strong solutions are provided. Several examples of Volterra equations
 having strong solutions are given, as well.
\end{abstract}

\section{Introduction} \label{AIs}

In this paper, which is the continuation of \cite{KaLi}, we
consider the following stochastic Volterra equation in a separable
Hilbert space $H$
\begin{equation} \label{AIe1}
 X(t) = X_0 + \int_0^t a(t-\tau) AX(\tau) d\tau  +
 \int_0^t \Psi(\tau)\,dW(\tau)\;.
\end{equation}
In (\ref{AIe1}), $X_0\in H$, $a\in
L^1_\mathrm{loc}(\mathbb{R}^+)$, $A$ is a closed unbounded linear
operator in $H$ with a dense domain $D(A)$ equipped with the graph
norm $|\cdot |_{D(A)}$. $W$ is a genuine Wiener process or a
cylindrical Wiener process and $\Psi$ is an appropriate process
defined below.

Equation (\ref{AIe1}) is motivated by a wide class of model
problems  and corresponds to an abstract stochastic version of
several deterministic problems, mentioned, e.g.\ in \cite{Pr} (see
also the references therein).

Let $(\Omega,\mathcal{F},(\mathcal{F}_t)_{t\geq 0},P)$ be a
stochastic basis and $U$ a separable Hilbert space. Let $Q\in
L(U)$ be a linear symmetric positive operator and $W(t),~t\geq 0$,
be an $U$-valued Wiener process with the covariance operator $Q$.
Let us note that the noise term $Z(t):=\int_0^t
\Psi(\tau)\,dW(\tau),~t\geq 0$, includes two distinguish cases.

When $\mathrm{Tr}\,Q<+\infty$, hence $W(t), ~t\geq 0$, is a
genuine Wiener process. Then we can take $U=H$, $\Psi:=I$ and the
noise term $Z(t)$ becomes $W(t),~t\geq 0$.

If $\mathrm{Tr}\,Q=+\infty, ~W(t), ~t\geq 0$, is so-called
cylindrical Wiener process. In this case, in order to provide a
sense of the integral $Z(t)$, the process $\Psi(t),~t\geq 0$, has
to be an operator-valued process (see, e.g.\ \cite{Ka}). We define
the subspace $U_0:=Q^{1/2}(U)$ of the space $U$ endowed with the
inner product $\langle u, v\rangle_{U_0}:=\langle Q^{-1/2}u,
Q^{-1/2}v\rangle_U$.

By $L_2^0:=L_2(U_0,H)$ we denote the set of all Hilbert-Schmidt
operators acting from $U_0$ into $H$; the set $L_2^0$ equipped
with the norm
$|C|_{L_2(U_0,H)}:=\left(\sum_{k=1}^{+\infty}|Cu_k|_H^2
\right)^{\frac{1}{2}}\!$, is a separable Hilbert space.

By $\mathcal{N}^2(0,T;L_2^0)$, where $T<+\infty$ is fixed, we
denote a Hilbert space of all
 $L_2^0$-predictable processes $\Psi$ such that
 $|| \Psi||_T < +\infty$, where
$$ 
||\Psi||_T := \left\{\mathbb{E}\left( \int_0^T
|\Psi(\tau)|_{L_2^0}^2\,d\tau \right) \right\}^{\frac{1}{2}} =
\left\{\mathbb{E} \int_0^T \left[ \mathrm{Tr}
(\Psi(\tau)Q^{\frac{1}{2}}) (\Psi(\tau)Q^{\frac{1}{2}})^* \right]
d\tau \right\}^{\frac{1}{2}}.
$$ 

If $\Psi\in \mathcal{N}^2(0,T;L_2^0)$, then the integral $\int_0^t
\Psi(\tau)\, dW(\tau)$ has sense.

In this paper, we use the so-called resolvent approach to the
Volterra equation (\ref{AIe1}) (for details we refer to
\cite{Pr}).

   \begin{definition}\label{AId1}
    A family $(S(t))_{t\geq 0}$ of bounded linear operators in
    a Banach space $B$ is
    called {\tt resol\-vent} for (\ref{AIe1}) if the following
    conditions are satisfied:
    \begin{enumerate}
     \item $S(t)$ is strongly continuous on $\mathbb{R}_+$ and $S(0)=I$;
     \item $S(t)$ commutes with the operator $A$: \\
     $S(t)(D(A))\subset D(A)$ and $~AS(t)x=S(t)Ax$
      for all $x\in D(A)$ and $t\geq 0$;
    \item the following {\tt resolvent equation} holds
    \begin{equation} \label{AIe2}
     S(t)x = x + \int_0^t a(t-\tau) A\,S(\tau)x \,d\tau \end{equation}
    for all $x\in D(A),~t\geq 0$.
    \end{enumerate}
   \end{definition}

Let us emphasize that the family $(S(t))_{t\geq 0}$ does not
create any semigroup and that $S(t), ~t\geq 0$, are generated by
the pair $(A,a(t))$, that is, the operator $A$ and the kernel
function $a(t), ~t\geq 0$.

A consequence of the strong continuity of $S(t)$ is that
$\sup_{t\leq T} \; ||S(t)|| < +\infty$ for any $T\geq 0$.

 \begin{definition}\label{AId2}
   We say that the function $a\in L^1(0,T)$ is {\tt completely positive}
   on $[0,T]$, if for any $\mu\geq 0$, the solutions of the equations
  \begin{equation} \label{AIe3}
   s(t) \!+\! \mu (a \star s)(t) =1 \quad \mbox{and} \quad r(t) \!+\! \mu (a
   \star r)(t) = a(t)
   \end{equation}
   satisfy $s(t)\geq 0$ and $ r(t) \geq 0$ on $[0,T]$.
  \end{definition}

 The class of
 completely positive kernels, introduced in \cite{ClNo}, arise naturally in
 applications, see~\cite{Pr}.

\begin{definition}\label{AId3}
Suppose $S(t),~t\geq 0$, is a resolvent.  $S(t)$ is called {\tt
exponentially bounded} if there are constants
 $M\geq 1$ and $\omega\in\mathbb{R}$ such that
 $$||S(t)|| \leq M\,e^{\omega t}, \mbox{~~for all~~} t\geq 0 \;;$$
 $(M,\omega)$ is called a {\tt type} of $S(t)$.
\end{definition}

Let us note that in contrary to $C_0$-semigroups, not every
resolvent family needs to be exponentially bounded; for
counterexamples we refer to \cite{DePr}.

In the paper, the key role is played by the following, yet
non-published, result providing a convergence of resolvents.

\begin{theorem}\label{AIt1}
   Let $A$ be the generator of a
   $C_0$-semigroup in $B$ and suppose
   the kernel function $a$ is completely po\-sitive. Then $(A,a)$
   admits an exponentially bounded resolvent $S(t)$. Moreover, there
   exist bounded operators $A_n$ such that $(A_n,a)$ admit resolvent
   families $S_n(t)$ satisfying $ ||S_n(t) || \leq Me^{w_0 t}~ (M\geq
   1,~w_0\geq 0)$ for all $t\geq 0,~n\in \mathbb{N},$ and
   \begin{equation}\label{AIe4}
     S_n(t)x \to S(t)x \quad \mbox{as} \quad n\to +\infty
   \end{equation}
   for all $x \in B,\; t\geq 0.$

Additionally, the convergence is uniform in $t$ on every compact
subset of $~\mathbb{R}_+$.
\end{theorem}

\begin{remark}

(a) The convergence (\ref{AIe4}) is an extension of the result due
to Cl\'ement \& Nohel \cite{ClNo}. The operators $A_n, ~n\in
\mathbb{N}$, are the Yosida approximation of the operator $A$. For
more details and the proof we refer to \cite{KaLi}.

(b) The above theorem give a partial answer to the following open
problem for a resolvent family $S(t)$ generated by a pair $(A,a)$:
do exist bounded linear operators $A_n$ generating resolvent
families $S_n(t)$ such that $S_n (t)x \to S(t)x$?. Note that in
case $a(t)\equiv 1$ the answer is yes, namely $A_n$ are provided
by the Hille-Yosida approximation of $A$ and $S_n(t) = e^{t A_n}.$
\end{remark}

\section{Probabilistic results} \label{AIs2}

In the sequel we shall use the following {\tt Probability
Assumptions}, abbr.\ (PA):
\begin{enumerate}
\item $X_0$ is an $H$-valued, $\mathcal{F}_0$-measurable random
variable;
\item $\Psi\in \mathcal{N}^2(0,T;L_2^0)$ and the interval
$[0,T]$ is fixed.
\end{enumerate}

The following types of definitions of solutions to (\ref{AIe1})
are possible, see \cite{Ka2}.

\begin{definition} \label{AId4}
Assume that (PA) hold. An $H$-valued predictable process
 $X(t),~t\in [0,T]$, is said to be a {\tt strong solution} to
(\ref{AIe1}), if $X$ takes values in $D(A)$, $P$-a.s.,
\begin{equation}\label{AIe5}
\int_0^T |a(T-\tau)AX(\tau)|_H \,d\tau<+\infty, \quad
P\mathrm{-a.s.}
\end{equation}
and for any $t\in [0,T]$ the equation (\ref{AIe1}) holds $P$-a.s.
\end{definition}

Let $A^*$ be the adjoint of $A$ with a dense domain
 $D(A^*)\subset H$ and the graph norm $|\cdot |_{D(A^*)}$.

\begin{definition} \label{AId5}
Let (PA) hold. An $H$-valued predictable process $~X(t),~t\in
[0,T]$, is said to be a {\tt weak solution} to (\ref{AIe1}), if
 $P(\int_0^t|a(t-\tau)X(\tau)|_H d\tau<+\infty)=1$ and if for all
 $\xi\in D(A^*)$ and all $t\in [0,T]$ the following equation holds
\begin{eqnarray*}
\langle X(t),\xi\rangle_H = \langle X_0,\xi\rangle_H + \langle
\int_0^t a(t-\tau)X(\tau)\,d\tau, A^*\xi\rangle_H \\ + \langle
\int_0^t \Psi(\tau)dW(\tau),\xi\rangle_H, \quad P\mathrm{-a.s.}
\end{eqnarray*}
\end{definition}

\begin{definition} \label{AId6}
Assume that $X_0$ is $\mathcal{F}_0$-measurable random variable
such that $P(X_0\in D(A))$ $=1$. An $H$-valued predictable process
 $X(t),~t\in [0,T]$, is said to be a {\tt mild solution} to the
stochastic Volterra equation (\ref{AIe1}), if $
\mathbb{E}(\int_0^t |S(t-\tau)\Psi(\tau)|_{L_2^0}^2 \,d\tau
)<+\infty$ for $t\leq T$ and, for arbitrary $t\in [0,T]$,
\begin{equation}\label{AIe6}
X(t) = S(t)X_0 + \int_0^t S(t-\tau)\Psi(\tau)\,dW(\tau), \quad
P\mathrm{-a.s.}
\end{equation}
\end{definition}

First, let us consider the case when $W$ is an $H$-valued genuine
Wiener process. In this case the equation (\ref{AIe1}) reads
\begin{equation}\label{AIe7}
X(t) = X_0 + \int_0^t a(t-\tau)\,AX(\tau)\,d\tau + W(t), \quad
t\geq 0\;.
\end{equation}
We define the convolutions:
\begin{eqnarray*}
 W_S(t) & := & \int_0^t S(t-\tau)\, dW(\tau)  \\
 W_{S_n}(t) & := & \int_0^t S_n(t-\tau)\, dW(\tau),
\end{eqnarray*}
where $S(t),~ S_n(t), ~t\geq 0$, are resolvents corresponding to
 $A$ and $A_n$, respectively.

 \medskip
 Now, we can recall several results from \cite{KaLi}, not published.

\begin{theorem}\label{AIt2}
Let $A$ be the generator of a $C_0$-semigroup in $H$. Suppose the
kernel function $a$ is completely positive and
 $\mathrm{Tr}\,Q<+\infty $. Then
$$
\lim_{n\rightarrow +\infty} \mathbb{E} ( \sup_{t\in [0,T]}
|W_S(t)-W_{S_n}(t)|_H^p )=0 $$ for any $p\geq 2$.
\end{theorem}

\begin{lemma} \label{AIl1}
Assume that $A$ is the generator of a $C_0$-semigroup in $H$, the
kernel function $a$ is completely positive and
 $\mathcal{R}(S(t))\subset D(A)$. If $X_0=0$, then the convolution
 $W_S(t)$ fulfills (\ref{AIe7}).
\end{lemma}

\begin{theorem} \label{AIt3}
Assume that $A$ is the generator of a $C_0$-semigroup and the
kernel function $a$ is completely positive. Let $
\mathcal{R}(S(t)) \subset D(A)$ for all $t>0$ and $X_0=0$. Then
the equation (\ref{AIe7}) has a strong solution. Precisely, the
convolution $W_S(t)$ is the strong solution to (\ref{AIe7}).
\end{theorem}

Theorem \ref{AIt3} is an extension of the semigroup case. It
provides sufficient conditions for mild solutions to be strong
solutions.

\medskip

Now, let us consider the case when $W$ is a cylindrical Wiener
process. We define the convolution
$$ W^\Psi(t) := \int_0^t S(t-\tau)\Psi(\tau)\,dW(\tau) $$
for $\Psi\in\mathcal{N}^2(0,T;L_2^0)$.

\begin{proposition} \label{AEp1}
If $\Psi\in\mathcal{N}^2(0,T;L_2^0)$ and
 $\Psi(\cdot,\cdot)(U_0)\subset D(A)$, $P$-a.s., then the
stochastic convolution $W^\Psi$ fulfills the equation
$$
 \langle W^\Psi(t),\xi\rangle_H =
 \int_0^t \langle a(t-\tau)W^\Psi(\tau), A^*\xi\rangle_H
 + \int_0^t \langle \xi,\Psi(\tau)dW(\tau)\rangle_H
$$
for any $t\in [0,T]$ and $\xi\in D(A^*)$.
\end{proposition}

Proposition \ref{AEp1} (see \cite{Ka2}) enables to formulate the
following results.
\begin{proposition} \label{AEp2}
Let $A$ be the generator of $C_0$-semigroup in $H$ and suppose
that the function $a$ is completely positive. If $\Psi$ and
 $A\Psi$ belong to $\mathcal{N}^2(0,T;L_2^0)$ and in addition
 $\Psi(\cdot,\cdot)(U_0)\subset D(A)$, $P$-a.s., then the following
equality holds
$$
W^\Psi(t) = \int_0^t a(t-\tau)A\, W^\Psi(\tau)\,d\tau
  + \int_0^t \Psi(\tau)\,dW(\tau)\;.$$
\end{proposition}

\begin{theorem} \label{AEt4}
Suppose that assumptions of Proposition 2 hold. Then the equation
(1) has a strong solution. Precisely, the convolution $W^\Psi$ is
the strong solution to~(\ref{AIe1}).
\end{theorem}

\section{Fractional Volterra equations} \label{AIs3}

Let us note that the condition $\mathcal{R}(S(t)) \subset D(A),
~t\geq 0$ used  in Theorem \ref{AIt3}, is satisfied by a large
class of resolvents. Particularly, when the equation (\ref{AIe7})
is parabolic in the sense of \cite{Pr} and $a(t)$ is $k$-regular,
e.g.\ $a(t)=t^{\alpha-1}/\Gamma(\alpha)$, $t\geq 0$, $\alpha\in
(0,2)$, where $\Gamma$ is the gamma function. This fact leads us
to fractional Volterra equation of the following form
\begin{equation}\label{AIe8}
 X(t) = X_0 + \int_0^t a(t-\tau) AX(\tau) d\tau  +
 \int_0^t \Psi(\tau)\,dW(\tau), \quad t\geq 0,
\end{equation}
when $a(t)=g_\alpha(t),\quad \alpha>0$ with
$g_\alpha(t)=\frac{t^{\alpha-1}}{\Gamma(\alpha)}$. Observe that
for $\alpha\in (0,1]$, $g_\alpha$ are completely positive, but for
$\alpha > 1$, $g_\alpha$ are not completely positive.

Now, the pairs $(A,g_\alpha(t))$ generate $\alpha$-times
resolvents $S_\alpha(t), ~t\geq 0$; for more details, see
\cite{Ba}.

\begin{remark} \label{rem1}
 Observe that the $\alpha$-times resolvent family corresponds to a
 $C_0$-semi\-group in case $\alpha=1$ and a cosine family in case $\alpha=2$.
 In consequence, when $1<\alpha<2$ such resolvent families interpolate
 $C_0$-semigroups and cosine functions. In particular, for $A=\Delta$,
 the integrodifferential equation corresponding to such resolvent family
 interpolates the heat equation and the wave equation, see, e.g.\ \cite{Fu}.
\end{remark}

We consider two cases:
\begin{description}
\item[(A1)]  $A$ is the generator of $C_0$-semigroup
 and $\alpha\in (0,1)$;
\item[(A2)]  $A$ is the generator of a strongly continuous
 cosine family and $\alpha\in (0,2)$.
\end{description}

In this part of the paper, the results concerning a weak convergence of
$\alpha$-times resolvents play the key role. Using the very recent result due to
Li and Zheng \cite{LiZh}, we can formulate the approximation theorems for
fractional Volterra equations.

\begin{theorem} \label{AIt5}
Let $A$ be the generator of a $C_0$-semigroup $(T(t))_{t\geq 0}$
in a Banach space $B$ such that
\begin{equation}\label{AIe9}
\|T(t) \| \leq Me^{\omega t}, \quad t \geq 0\;.
\end{equation}
Then, for each $ 0 < \alpha < 1$
there exist bounded operators $A_n$ and $\alpha$-times resolvent
families $S_{\alpha,n}(t)$ for $A_n$ satisfying $
||S_{\alpha,n}(t) || \leq MCe^{(2\omega)^{1/\alpha} t}$, for all
$t\geq 0,~n\in \mathbb{N}$, and
\begin{equation} \label{AIe10}
S_{\alpha,n}(t)x \to S_{\alpha}(t)x \quad \mbox{as} \quad n\to
+\infty
\end{equation}
for all $x \in B,\; t\geq 0$. Moreover, the convergence is uniform
in $t$ on every compact subset of $ \mathbb{R}_+$.
\end{theorem}

Outline of the proof: The first assertion follows  from
\cite[Theorem 3.1]{Ba}, that is, for each $ 0<\alpha < 1$ there is
an $\alpha$-times resolvent family $(S_{\alpha}(t))_{t\geq 0}$ for
$A$ given by
$$ 
S_{\alpha}(t)x = \int_0^{\infty} \varphi_{t,\alpha}(s) T(s)xds,
\quad t>0,
$$ 
where $\varphi_{t,\gamma}(s) := t^{-\gamma}
\Phi_{\gamma}(st^{-\gamma})$ and $ \Phi_{\gamma}(z)$ is the Wright
function defined as
$$ 
\Phi_{\gamma}(z):= \sum_{n=0}^{\infty} \frac{(-z)^n}
{n!\,\Gamma(-\gamma n + 1 - \gamma)}, \quad 0 < \gamma < 1.
$$ 
Define
$$ 
A_n := n AR(n,A) = n^2 R(n,A) - nI, \qquad n> w,
$$ 
the {\it Yosida approximation} of $A$.

Since each $A_n$ is bounded, it follows that for each $0< \alpha <
1$ there exists an $\alpha$-times resolvent family $
(S_{\alpha,n}(t))_{t\geq 0}$ for $A_n$ given as

$$ 
S_{\alpha,n}(t) = \int_0^{\infty} \varphi_{t,\alpha}(s)
e^{sA_n}ds, \quad t>0.
$$ 

We recall that the Laplace transform of the Wright function
corresponds to $E_{\gamma}(-z)$ where $E_{\gamma}$ denotes the
Mittag-Leffler function. In particular, $\Phi_{\gamma}(z)$ is a
probability density function. It follows that for $t\geq 0$:
\begin{eqnarray*}
\| S_{\alpha,n}(t) \| &\leq & \int_0^{\infty}
\varphi_{t,\alpha}(s) \| e^{s A_n} \| ds \\ &\leq& M
\int_0^{\infty} \varphi_{t,\alpha}(s) e^{2\omega s} ds = M
\int_0^{\infty} \Phi_{\alpha}(\tau) e^{2\omega t^{\alpha} \tau }
d\tau= M E_{\alpha}(2 \omega t^{\alpha}).
\end{eqnarray*}

The continuity in $t\geq 0$ of the Mittag-Leffler function
and its asymptotic behavior, imply that for
$\omega \geq 0$ there exists a constant $C>0$ such that
$$ 
E_{\alpha}(\omega t^{\alpha}) \leq C e^{{\omega^{1/\alpha}} t},
\quad t \geq 0,\,\, \alpha \in (0,2).
$$ 

This gives

$$ 
\| S_{\alpha,n}(t) \| \leq MCe^{(2\omega)^{1/\alpha}t}, \quad t
\geq 0.
$$ 

Now we recall the fact that $ R(\lambda,A_n)x \to R(\lambda,A)x $
as $ n\to \infty$ for all $\lambda $ sufficiently large (e.g.\
\cite[Lemma~7.3]{Pa}), so we can conclude from \cite[Theorem
4.2]{LiZh} that
$$ 
S_{\alpha,n}(t)x \to S_{\alpha}(t)x \quad \mbox{as} \quad n\to
+\infty
$$ 
for all $x \in B,$ uniformly for $t$ on every compact subset of
$\mathbb{R}_+$ \hfill $\square$ \\

An analogous convergence for $\alpha$-times resolvents can be proved in another
case, too.

\begin{theorem} \label{AIt6}
Let $A$ be the generator of a $C_0$-cosine family $(T(t))_{t\geq
0}$ in a Banach space $B$. Then, for each $0<\alpha<2$
there exist bounded operators $A_n$ and $\alpha$-times resolvent
families $S_{\alpha,n}(t)$ for $A_n$ satisfying $
||S_{\alpha,n}(t) || \leq MCe^{(2\omega)^{1/\alpha} t}$, for all
$t\geq 0,~n\in \mathbb{N}$, and $ S_{\alpha,n}(t)x \to
S_{\alpha}(t)x$ as $~n\to +\infty$ for all $x \in B,\; t\geq 0$.
Moreover, the convergence is uniform in $t$ on every compact
subset of $ \mathbb{R}_+$.
\end{theorem}

Now, we are able to formulate the results analogous to that in
section \ref{AIs2}.

\begin{theorem} \label{AIt7}
Assume that {\bf (A1)} or {\bf (A2)} holds and
$\mathcal{R}(S_\alpha(t))\subset D(A)$. If $X_0=0$, then the equation
(\ref{AIe7}) has a strong solution.
\end{theorem}

Outline of the proof:
We define the convolutions
\begin{eqnarray} \label{eSW12}
W_{S_\alpha}(t) & := & \int_0^t S_\alpha(t-\tau)\, dW(\tau), \\
 \label{eSW13}
W_{S_{\alpha,n}}(t) & := & \int_0^t S_{\alpha,n}(t-\tau)\,
dW(\tau),
\end{eqnarray}
where $S_\alpha(t), S_{\alpha,n}(t),~t\geq 0$, are resolvents
corresponding to the operators $A$ and $A_n$, respectively.

By \cite[Corollary 1]{Ka2}, for every $n\in\mathbb{N}$ the
convolution $W_{S_{\alpha,n}}(t)$ fulfills (\ref{AIe7}), where
$W_{S_{\alpha,n}}(t)$ and $A_n$ are as above. Because
Theorem \ref{AIt2} holds for the convolutions $W_{S_{\alpha,n}}(t)$
and  $W_{S_{\alpha}}(t)$,
we can assume that
$W_{S_{\alpha,n}}(t)\to W_{S_\alpha}(t)$ as $n\to +\infty$, P-a.s.

By the Lebesgue dominated convergence theorem
\begin{equation}\label{eSW16a}
\lim_{n\to +\infty} \; \sup_{t\in [0,T]} \; \mathbb{E}
|W_{S_{\alpha,n}}(t)-W_{S_\alpha}(t)|_H^2 =0.
\end{equation}
From assumptions $\mathcal{R}(S_\alpha(t))\subset D(A)$, so
$P(W_{S_\alpha}(t)\in D(A))=1$.

Using (\ref{eSW16a}) and the fact that $\displaystyle
\lim_{n\to\infty} A_n x =Ax$ for any $x\in D(A)$, we have
$$ \lim_{n\to +\infty} \; \sup_{t\in [0,T]} \; \mathbb{E}
|A_n W_{S_{\alpha,n}}(t)-A W_{S_\alpha}(t)|_H^2 =0. $$ Hence,
 $W_{S_\alpha}(t)=\int_0^t
g_\alpha(t-\tau)\,AW_{S_\alpha}(\tau) \,d\tau + W(t), ~ t\in
[0,T]$, P-a.s.

Because the convolution $W_{S_\alpha}(t)$ has integrable trajectories (see
\cite{Ka2}) and  the
closed linear unbounded operator $A$ becomes bounded on $D(A)$
endowned with the norm $|\cdot |_{D(A)}$ (see \cite[Chapter
5]{We}), then $A\,W_{S_\alpha}(\cdot )\in L^1([0,T];H)$. Hence,
the function $g_\alpha(T-\tau)\,A\,W_{S_\alpha}(\tau)$
is integrable with respect to $\tau$, what finishes the proof.
\hfill $\square$

\begin{theorem} \label{AIt8}
Assume that {\bf (A1)} or {\bf (A2)} holds. If $\Psi$ and $A\Psi$
belong to $\mathcal{N}^2(0,T;L_2^0)$ and in addition
$\Psi(\cdot,\cdot)(U_0)\subset D(A)$, $P$-a.s., then the equation
(\ref{AIe1}) with $X_0=0$ has a strong solution. Precisely, the
convolution
$$W_\alpha^\Psi(t):=\int_0^t S_\alpha (t-\tau)\,\Psi(\tau)\,dW(\tau)$$
 is the strong solution to (\ref{AIe1}).
\end{theorem}

Outline of the proof:
First, analogously like in \cite{Ka2}, we show that the convolution
$W_\alpha^\Psi(t)$ fulfills the following equation
\begin{equation} \label{eSW19}
W_\alpha^\Psi(t) = \int_0^t g_\alpha(t-\tau)A\,
W_\alpha^\Psi(\tau)\,d\tau
  + \int_0^t \Psi(\tau)\,dW(\tau)\,.
\end{equation}
Next, we have to show the condition
\begin{equation} \label{eSW3.1}
\int_0^T |g_\alpha(T-\tau)AW_\alpha^\Psi(\tau)|_H \,d\tau<+\infty,\quad
P-a.s., \quad \alpha >0,
\end{equation}
that is, the condition (\ref{AIe5}) adapted for the fractional Volterra
equation (\ref{AIe8}).

The convolution $W_\alpha^\Psi (t)$ has integrable trajectories (see
\cite{Ka2}), that is, $W_\alpha^\Psi (\cdot )\in L^1([0,T];H)$,
P-a.s. The closed linear unbounded operator $A$ becomes bounded on
($D(A),|\cdot|_{D(A)}$), see \cite[Chapter 5]{We}. Hence,
$AW_\alpha^\Psi (\cdot )\in L^1([0,T];H)$, P-a.s. Therefore, the
function $g_\alpha(T-\tau)AW_\alpha^\Psi (\tau)$ is integrable
with respect to $\tau$, what completes the proof.
\hfill $\square$

\section{Examples}\label{AIs4}

The class of equations fulfilling our conditions depends on where
the operator $A$ is defined, in particular, the domain of $A$
depends on each considered problem, and also depends on the
properties of the kernel function $a$.

Let $\Omega$ be a bounded domain in $\mathbb{R}^n$ with smooth
boundary $\partial \Omega$. Consider the differential operator of
order $2m$:
\begin{equation}{\label{eq1}}
A(x,D)= \sum_{|\alpha|\leq 2m } a_{\alpha}(x) D^{\alpha}
\end{equation}
where the coefficients $a_{\alpha}(x)$ are sufficiently smooth
complex-valued functions of $x$ in $\overline \Omega.$ The
operator $A(x,D)$ is called {\tt strongly elliptic} if there
exists a constant $c>0$ such that
\begin{equation}
Re(-1)^m \sum_{|\alpha|= 2m } a_{\alpha}(x) \xi^{\alpha} \geq c
|\xi|^{2m}
\end{equation}
for all $ x\in \overline \Omega$ and $\xi \in \mathbb{R}^n.$

Let $A(x,D)$ be a given strongly elliptic operator on a bounded
domain $\Omega \subset \mathbb{R}^n$ and set $D(A)= H^{2m}(\Omega)
\cap H_0^{m}(\Omega).$ For every $u \in D(A)$ define
\begin{equation}
Au = A(x,D)u.
\end{equation}
Then the operator $-A$ is the infinitesimal generator of an
analytic semigroup of operators on $H= L^2(\Omega)$ (cf.
\cite[Theorem 7.2.7]{Pa}). We note that if the operator $A$ has
constant coefficients, the result remains true for the domain
$\Omega = \mathbb{R}^n.$

A concrete example is the Laplacian
\begin{equation}
\Delta u = \sum_{i=1}^n \frac{\partial^2 u}{\partial x_i^2},
\end{equation}
since $- \Delta$ is clearly strongly elliptic. It follows that
$\Delta u$ on $D(A) = H^2( \Omega) \cap H_0^1(\Omega)$ is the
infinitesimal generator of an analytic semigroup on $L^2(\Omega)$.

In particular, by \cite[Corollary 2.4]{Pr} the operator $A$
given by (\ref{eq1}) generates an analytic resolvent $S(t)$
whenever $a \in C(0,\infty) \cap L^1(0,1)$ is completely
monotonic. As a consequence $\mathcal{R}(S(t)) \subset D(A)$ for
all $t>0.$

This example fits in our results (Theorem \ref{AIt3}) if $a$ is
also completely positive. For example: $a(t) =
t^{\alpha-1}/\Gamma(\alpha)$ is both, completely positive and
completely monotonic for $0 < \alpha \leq 1$ (but not for $\alpha
>1$).

Another  class of examples is provided by the following: suppose
$a \in L^1_{loc}(\mathbb{R}_+ )$ is of subexponential growth and
$\pi/2$-sectorial, and let $A$ generate a bounded analytic
$C_0$-semigroup  in a complex Hilbert space $H$.  Then it follows
from \cite[Corollary 3.1]{Pr} that the Volterra equation of scalar
type $ u= a*Au +f$ is parabolic. If, in addition, $a(t)$ is
$k$-regular for all $k\geq 1$ we obtain from \cite[Theorem 3.1
]{Pr} the existence of a
 resolvent $S\in C^{k-1}((0,\infty),
\mathcal{B}(H))$ such that $ \mathcal{R}(S(t)) \subset D(A)$ for
$t>0$ (see \cite[ p.82 (f)]{Pr}). These observations  together with
Theorem \ref{AIt3}  give us the following result.

\begin{corollary}
Suppose that $A$ generates a bounded analytic $C_0$-semigroup in
a complex Hilbert space $H$ and $a \in L^1_{loc}(\mathbb{R}_+ )$
is of subexponential growth, $\pi/2$-sectorial, completely
positive and $k$-regular for all $k\geq 1.$ Then the equation
$$ 
X(t) = X_0 + \int_0^t a(t-\tau)\,AX(\tau)\,d\tau + W(t), \quad
t\geq 0\;.
$$ 
 has a strong solution.
\end{corollary}

\medskip

 \end{document}